\documentclass[12pt]{amsart}
\usepackage{amsmath,amssymb,amsfonts,amsthm,amstext,verbatim}
\usepackage{url,enumerate,color,graphicx,stmaryrd,psfrag} 
\usepackage{hyperref}
\usepackage{ctable,multirow}
\hypersetup{colorlinks=true, linkcolor=darkred, citecolor=darkgreen, urlcolor=darkblue}
\definecolor{darkred}{RGB}{150,0,0}
\definecolor{darkgreen}{RGB}{0,150,0}
\definecolor{darkblue}{RGB}{0,0,200}
\newcommand\noi{\noindent}
\newcommand\resp{respectively}
\newcommand\seprule{\noalign{\vskip0.5ex}\hline\noalign{\vskip0.5ex}}
\newcommand\bp{base$+$prop}
\newcounter{mycount}

\newenvironment{numlist}{\begin{list}{\arabic{mycount}.}%
   {\usecounter{mycount}\labelwidth=1cm\itemsep 0pt}}{\end{list}}

\newenvironment{dashlist}{\begin{list}{$\circ$}%
   {\usecounter{mycount}\labelwidth=1cm\itemsep 0pt}}{\end{list}}
\newenvironment{Letlist}{\begin{list}{\rm\Alph{mycount}.}%
   {\usecounter{mycount}\labelwidth=1cm\itemsep 0pt}}{\end{list}}

\begin{document}
\title[European Apportionment via the Cambridge Compromise]{European Apportionment via the{\\}Cambridge Compromise}
\author{Geoffrey R.\ Grimmett}
\address{Statistical Laboratory,\\Centre for
Mathematical Sciences, Cambridge University,\\Wilberforce Road,
Cambridge CB3 0WB, UK} 
\email{g.r.grimmett@statslab.cam.ac.uk}
\urladdr{http://www.statslab.cam.ac.uk/$\sim$grg/}

\begin{abstract}
Seven mathematicians and one political scientist met at the Cambridge Apportionment Meeting 
in January 2011.
They agreed a unanimous recommendation to the
 European Parliament for its future apportionments between the
EU Member States. This is a short factual account of the reasons
that led to the Meeting, of its debates and
report, and of some of the ensuing Parliamentary debate.
\end{abstract}

\date{Revised August 20,  2011}
\renewcommand{\subjclassname}{\textup{2010} Mathematics Subject
Classification}

\keywords{Apportionment problem, European Parliament, degressive proportionality,
\bp\ method, D'Hondt method}
\subjclass[2010]{91B12}

\maketitle
\section{Background and Brief}

\subsection{Background}
As the European Union has grown and its population has developed, so has the 
constitution and structure of the European Parliament.  In recognition of the need for an orderly
allocation of Parliamentary seats between the EU Member States, 
its  Committee on Contitutional Affairs (AFCO)
commissioned a Symposium of Mathematicians
to \lq\lq identify a mathematical formula
for the distribution of seats which will be durable, transparent and impartial to politics''. 
The purposes of the reform were described thus in \cite{duff}:  
\begin{dashlist}
\item
The aim of the symposium is to discuss and, if possible, to propose to the Committee on
Constitutional Affairs a mathematical formula for the redistribution of the 751 seats in
the European Parliament. The formula should be as transparent as possible and capable
of being sustained from one Parliamentary mandate to the next.
\item
The purpose of the Symposium is to eliminate the political bartering which has characterised
the distribution of seats so far by enabling a smooth reallocation of seats once every
five years which takes account of migration, demographic shifts and the accession of
new Member States.
\end{dashlist}

The current note is more a record of the events surrounding the Cambridge Apportionment Meeting than
it is a critical analysis of the politics.
An account of the history of the current apportionment of Parliament, and of the associated
\lq\lq political bartering",  may be found in \cite{duff2}.

\subsection{Cambridge Apportionment Meeting (CAM)}
The Symposium took place in the Centre for Mathematical Sciences,
Cambridge University, on 28--29 January 2011. The following participated:
\begin{dashlist}
\item Geoffrey Grimmett (University of Cambridge), Director,
\item Friedrich Pukelsheim (University of Augsburg), co-Director,
\item Jean-Fran\c cois Laslier (\'Ecole Polytechnique, Paris),
\item Victoriano Ram\'\i rez Gonz\'alez (University of Granada),
\item Richard Rose (University of Aberdeen; European University Institute, Florence),
\item Wojciech S\l omczy\'nski (Jagiellonian University, Krak\'ow),
\item Martin Zachariasen (University of Copenhagen),
\item Karol \.{Z}yczkowski (Jagiellonian University, Krak\'ow),
\end{dashlist}
\noi
advised by
\begin{dashlist}%AFCO Committee in attendance
\item Andrew Duff MEP (AFCO Rapporteur),
\item Rafa\l\ Trzaskowski MEP (AFCO Vice-President),
\item Guy Deregnaucourt (AFCO Administrator),
\item Wolfgang Leonhardt (AFCO Administrator),
\item Kevin Wilkins (Assistant to Andrew Duff),
\end{dashlist}
\noi
in the presence of
\begin{dashlist}
\item Thomas Kellermann (College of Europe, Natolin, Warsaw),
\item Kai-Friederike Oelbermann (University of Augsburg).
\end{dashlist}

The formal Report of the Cambridge Apportionment Meeting
to the Congressional Affairs Committee  may be found  at \cite{camcom}.
The discussions and recommendations of CAM are summarized in the current article, together
with an account of some of the subsequent debate within the Committee.
Opinions expressed here are those of the author alone.

\subsection{The constraints}\label{sec:const}
Seat allocations are currently required to adhere to the terms of the Treaty of Lisbon.
\begin{dashlist}
\item Each Member State is to receive a minimum of 6 seats,
\item  and a maximum of 96 seats,
\item Parliament is constrained to have no more than 751 seats in total (including 
that of the President),
\item allocations are required to satisfy a condition of \lq\lq degressive proportionality".
\end{dashlist}

CAM was advised by the 
AFCO representatives that the first three constraints
are indeed \emph{inequalities} rather than \emph{equalities}, 
but nevertheless there existed a general expectation in Parliament
that its total size should not be less than 751, and that the smallest States should receive
an allocation not greater than 6 seats.
The issue of \lq\lq degressive proportionality" is formulated in more detail in Section \ref{sec:dp0}.
In reaching its conclusions, the Symposium took into account the following additional
observations concerning the general structure of the European Parliament:
\begin{dashlist}
\item
the EU has currently 27 Member States,
\item the smallest population (as published officially by 
Eurostat\footnote{\url{http://epp.eurostat.ec.europa.eu/}} is currently 412,970,
and the largest 81,802,257,
\item future accessions may include a number of States with a spread of sizes,
\item there will be migration and demographic changes,
\item Member States' population figures will be used as input
to the formula.
\end{dashlist}

\subsection{The criteria}\label{sec:crit}

Participants were sensitive in discussions to the three descriptors
 provided by the AFCO Committee, namely that the \lq \lq formula'' was required to be \emph{durable}, 
\emph{transparent}
and \emph{impartial to politics}.

\noi
\emph{Durable}:  A formula that adapts naturally to possible structural changes
in the architecture of the EU, arising for example through accessions by new States,
through migration, or through demographic shifts.

\noi
\emph{Transparent}: An apportionment method that is capable of simple and reasonable explanation
to EU citizens, irrespective of their backgrounds.

\noi
\emph{Impartial to politics}: A principled and fresh approach, unprejudiced with respect to particular
Member States or Political Groups, and free of influence 
from historical positions beyond the constraints of
Section \ref{sec:const}.

\subsection{Summary}

A discussion of degressive proportionality is to be found in
Section \ref{sec:dp0}. Section \ref{sec:cc}  contains a discussion of the main
recommendations of the Cambridge Apportionment Meeting, which are listed
explicitly 
in Section \ref{sec:recom}. A brief account of the subsequent
debate and resolutions of the Committee on Constitutional Affairs is
presented in Section \ref{sec:afco}. This chapter in the story of European
Apportionment ends with the shelving of the mathematical approach.

\section{Degressive Proportionality}\label{sec:dp0}

\subsection{Lamassoure--Severin definition}\label{sec:dp}

\emph{Degressive proportionality} has been defined in Paragraph 6 of the 
Lamassoure--Severin (2007) Motion of \cite{lamsev} as follows.
\begin{numlist}
\addtocounter{mycount}{5}
\item
  {[The European Parliament]} \lq\lq Considers that the principle of degressive
proportionality means that the ratio between the population and the number of
seats of each Member State must vary in relation to their respective populations
in such a way that each Member from a more populous Member State represents
more citizens than each Member from a less populous Member State and
conversely, but also that no less populous Member State has more seats than a
more populous Member State.''
\end{numlist}

The principle of degressive proportionality attracted significant debate and
a  major recommendation at CAM. 

\subsection{CAM recommendation}\label{sec:dp2}

It was noted that \emph{degressive proportionality} comprises two requirements:
\begin{numlist}
\item no smaller State shall receive more seats than a larger State,
\item the ratio population/seats shall increase as population increases.
\end{numlist}
Condition 1 is easy to accept.  Condition 2, on the other hand,  poses a
serious practical difficulty, and has in addition been violated in recent Parliamentary
apportionments. As noted in \cite{ma-rg08,puk10,RG10,rg-pal-m06}
and elsewhere, there are hypothetical instances of apportionment
for which there exists no solution satisfying both Condition 1 and Condition 2. 
There was an extensive discussion of this issue at CAM, 
centred on the following two Options.
\begin{Letlist}
\item Adopt a method whose outcomes invariably
satisfy Condition 2 but with a possibly reduced
Parliament-size.
\item Propose a change to the Lamassoure--Severin definition of degressive proportionality
 lying within existing law and allowing greater flexibility and transparency.
\end{Letlist}
A method satisfying Option A was presented at CAM (and is summarized in \cite[Sect.\ 6.2]{camcom}).
However, CAM preferred Option B on the
grounds of transparency of method, and the desirability
of achieving a given Parliament-size.

The recommendation of CAM was to amend Paragraph 6 of the Lamassoure--Severin Motion \cite{lamsev}
through the addition of the italicized phrase as follows.
\begin{numlist}
\addtocounter{mycount}{5}
\item
{[The European Parliament]} Considers that the principle of degressive proportionality
means that the ratio between the population and the number of seats of each
Member State \emph{before rounding to whole numbers} must vary in relation to their
respective populations in such a way that each Member from a more populous
Member State represents more citizens than each Member from a less populous
Member State and conversely, but also that no less populous Member State has
more seats than a more populous Member State.
\end{numlist}

\section{Cambridge Compromise}
\label{sec:cc}

\subsection{Base+prop method}\label{sec:ccbp}

The `Cambridge Compromise' recommendation\footnote{The Cambridge Compromise
proposal is named in harmony with the so-called Jagiellonian 
Compromise proposal of \cite{slomz10,slomz04} for voting within the European Commission.}
to the European Parliament is to adopt a \bp\ 
system, formulated in \cite{puk10}  as follows. 

The \emph{\bp} method proceeds in two stages. At the first stage, a fixed \emph{base}
number of seats is allocated to each Member State. At the second stage, the remaining
seats are allocated to States in \emph{prop}ortion to their population-sizes (subject to rounding,
and capping at the maximum). In order to achieve the given Parliament-size, one
introduces a further ingredient called the \emph{divisor}.

For given \emph{base} $b$, \emph{maximum} $M$, and \emph{divisor} $d$, 
define the associated \emph{allocation function} $A_d: [0,\infty) \to [0,\infty)$ by
$$
A_d(p) = \min\bigl\{b+p/d,M\bigr\},
$$
The \bp\ method is formulated as follows in mathematical terms. 
\begin{numlist}
\item
Assign to a Member State with population $p$ the \emph{seat share} $A_d(p)$,
\item
perform a rounding of the seat share $A_d(p)$ into an integer \emph{seat number}
$[A_d(p)]$,
\item
adjust the divisor $d$ in such a way that the sum of the seat numbers of all Member
States equals the given Parliament-size.
\end{numlist}

The total house-size with divisor $d$ is
$$
T(d) = \sum_i [A_d(p_i)],
$$
where the summation is over all Member States. 
The value of $d$ is chosen in such a way that $T(d)$ 
equals the prescribed total\footnote{There is
normally an interval of such $d$-values, and there are standard approaches to the question
of so-called \emph{ties}. See \cite{baly}, for example, and also Section \ref{sec:dhondt}.}.

The CAM recommendation is to use the base $b= 5$, and to use \emph{rounding
upwards}. Outcomes of the Cambridge Compromise
are presented in Tables \ref{table1} and \ref{table2}, 
with 2011 population figures taken
from the Eurostat website, and with 27, 28, and 29 Member States. 

It was through principled discussion that this recommendation was reached;
CAM was instructed to overlook historical apportionments, including the \emph{status quo}
as listed in Table \ref{table2}. Participants recognised the challenges that can be presented by change,
and these challenges proved formidable for the AFCO Committee (see Section \ref{sec:afco}).  

\subsection{Why \bp?}

The CAM participants considered a variety of 
apportionment schemes based around several different
linear and non-linear apportionment 
functions\footnote{Note that every non-decreasing concave apportionment
function leads invariably to allocations satisfying the revised form of degressive 
proportionality of Section \ref{sec:dp2}.}.
Linear functions were preferred over non-linear functions on grounds of transparency
and greater potential for proportionality.  The dual constraints of maximum and house-size 
are obstacles to the search for a \emph{smooth} linear apportionment function (that
is, a function
that is continuously differentiable, say).

Non-linear apportionment functions (following a power\footnote{A power-weighted
variant of the Cambridge Compromise is analysed in \cite{gop}.} or parabolic law, for example) can 
accommodate numerical constraints in a smoother manner.  
They can be used to fit curves to plots of data points distributed 
along (possibly concave) 
lines of trends, such as the current allocations to Member States. 
On the other hand,  they suffer from arbitrariness, 
and from lack of transparency.  The exercise confronting CAM was not
one of fitting a curve to historic data, but rather to form a fresh view of apportionment
that is impartial to yesterday's politics. 

From amongst linear systems, the \emph{\bp} method  stands out for its transparency.
It is degressively proportional in an active way, since
the base operates to the profit  of Member States at the lower end of the population table.
CAM considered that its simplicity outweighed the discontinuity in the
first derivative that arises currently through the maximum cap of 96 seats. We noted that
this discontinuity will diminish as the EU changes its shape through accessions. 
The recommendation to adopt the \bp\ method was reached through consideration 
of durability, transparency, impartiality, and degressive proportionality.

CAM noted in passing
that the \bp\ method can be interpreted as one in which the base is an 
allocation to Member States, and the remaining seats (prop) are 
proportional to population (subject to capping at the maximum). 
This resonates with the founding principles of the EU, enshrined in the Treaty, that the 
Union is made up both of Member States (enjoying equality in international law) 
and of citizens (enjoying democratic equality).

\subsection{Choice of base and rounding method}\label{sec:choice}

The choices of base and rounding methods are intertwined. A 
smaller base tends to favour larger States;
rounding \emph{upwards} is usually viewed as tending to favour 
smaller States.
These choices are informed by the existence of a minimum number 
$m$ of seats per State, and by degressive proportionality. 

Let us write $b+$R to denote the system with base $b$ and rounding method R, 
where R may denote one of:
\begin{itemize}
\item [{U:}] upwards rounding, 
\item [{S:}] standard rounding to the nearest integer,
\item [{D:}] downwards rounding. 
\end{itemize}
We say that the roundings of a real number $x$
are \emph{well defined} if $x$ is not an integer multiple of $\frac12$.
It was considered preferable, in the interests of transparency,
that the base be an integer.

Recall that $m=6$, and there is an expectation
that the smallest States will indeed receive 6 seats. 
It was therefore natural
to concentrate on the two possibilities: 
\begin{list}{}{\setlength{\labelwidth}{1cm}  \setlength{\leftmargin}{2cm}}
\item [{6+S:}] base $b=6$, \emph{standard rounding} (S),
\item [{5+U:}] base $b=5$, \emph{upwards rounding} (U).
\end{list}
Each of these two systems allocates at least 6 seats to every State. The minimum allocation
is however fragile under the first system (6+S), as 
illustrated in \cite[Sect.\ 5.3]{camcom} as follows. 
The currently smallest Member State is Malta, with a population of 412,970, and it receives an allocation
of 6 seats under both the above systems. If, however, its population were to increase by only 8,000
(other populations remaining unchanged), its allocation under 6+S rises to 7. This was considered unacceptable,
and for this reason CAM recommended 5+U.

There is an explicit trade-off between base and rounding method 
(see \cite{sj11,sj11b,zach}). Let $x$ be a real
number, and let $\lfloor\cdot\rfloor$ (\resp,  $\lceil\cdot\rceil$, $[\cdot]$) denote
 rounding downwards (\resp, upwards, and to the nearest integer). For any `base' $b$,
we have
$$
\lceil b+x\rceil = [b+\tfrac12+x] = \lfloor b+1+x\rfloor,
$$
whenever the roundings are well defined.
Subject to the last assumption, 
the three systems $b+$U, $(b+\frac12)+$S, $(b+1)+$D
result in the same allocations. 
In this sense, the systems 5+U, $5\frac12$+S, 6+D are equivalent.

\begin{table}[htbp]

\centering
\begin{tabular*}{\textwidth}{@{\extracolsep{\fill}} cc r c @{}r r}
\toprule
&\multirow{3}{*}{\it Member State}&\multirow{3}{*}{\it Population}&\multirow{3}{*}{\it Seats}&\multicolumn{1}{c}{\it Popn/seats}&\multicolumn{1}{c}{\it Popn/seats}\\
&&&&\multicolumn{1}{c}{\it before}&\multicolumn{1}{c}{\it after}\\
&&&&\multicolumn{1}{c}{\it rounding}&\multicolumn{1}{c}{\it rounding}\\
\midrule
1&Germany&81,802,257&96&  852,106.8 &852,106.8\\
2&France&64,714,074&85&770,259.3&\emph{761,342.0}\\
3&UK&62,008,048&81&768,264.0&765,531.5\\
4&Italy&60,340,328&79&766,950.8&763,801.6\\
5&Spain&45,989,016&62&752,036.4&741,758.3\\
\seprule
6&Poland&38,167,329&52&739,643.2&733,987.1\\
7&Romania&21,462,186&32&687,772.5&670,693.3\\
8&Netherlands&16,574,989&26&656,745.2&637,499.6\\
9&Greece&11,305,118&19&601,222.1&595,006.2\\
10&Belgium&10,839,905&19&594,438.5&\emph{570,521.3}\\
\seprule
11&Portugal&10,637,713&18&591,356.6&590,984.1\\
12&Czech Rep.&10,506,813&18&589,315.9&583,711.8\\
13&Hungary&10,014,324&18&581,298.7&556,351.3\\
14&Sweden&9,340,682&17&569,380.7&549,451.9\\
15&Austria&8,375,290&16&550,056.4&523,455.6\\
\seprule
16&Bulgaria&7,563,710&15&531,334.8&504,247.3\\
17&Denmark&5,534,738&12&470,724.2&461,228.2\\
18&Slovakia&5,424,925&12&466,706.8&452,077.1\\
19&Finland&5,351,427&12&463,965.8&445,952.2\\
20&Ireland&4,467,854&11&427,330.9&406,168.5\\
\seprule
21&Lithuania&3,329,039&10&367,250.6&332,903.9\\
22&Latvia&2,248,374&8&290,290.0&281,046.8\\
23&Slovenia&2,046,976&8&272,953.4&255,872.0\\
24&Estonia&1,340,127&7&201,939.0&191,446.7\\
25&Cyprus&803,147&6&134,291.1&133,857.8\\
\seprule
26&Luxembourg&502,066&6&89,446.6&83,677.7\\
27&Malta&412,970&6&75,027.7&68,828.3\\ [1ex]
\midrule
&{\it Total}&501,103,425&751& &\\
\bottomrule
\end{tabular*}
\bigskip
\small\caption{Each State receives one non-base seat for every 819,000 citizens or part thereof.
Population/seat ratios are strictly decreasing before rounding, but there are two violations after rounding,
namely Belgium and France when reading for the bottom. Data in this and the next table are
taken from the Eurostat website \url{http://epp.eurostat.ec.europa.eu/}.}
\label{table1}
\end{table}

\begin{table}[htbp]

\centering
\begin{tabular*}{\textwidth}{@{\extracolsep{\fill}} cc r c c c c}
\toprule
&\multirow{2}{*}{\it Member State}&\multirow{2}{*}{\it Population}&\multirow{2}{*}{\it Now}&{\it Seats}&{\it Seats}&{\it Seats}\\
&&&&{\it 27 States}&{\it 28 States}&{\it 29 States}\\
\midrule
1&Germany&81,802,257&99&96&96&96\\
2&France&64,714,074&74&85&83&82\\
3&UK&62,008,048&73&81&80&79\\
4&Italy&60,340,328&73&79&78&77\\
5&Spain&45,989,016&54&62&61&60\\
\seprule
6&Poland&38,167,329&51&52&51&51\\
7&Romania&21,462,186&33&32&31&31\\
8&Netherlands&16,574,989&26&26&25&25\\
9&Greece&11,305,118&22&19&19&19\\
10&Belgium&10,839,905&22&19&18&18\\
\seprule
11&Portugal&10,637,713&22&18&18&18\\
12&Czech Rep.&10,506,813&22&18&18&18\\
13&Hungary&10,014,324&22&18&17&17\\
14&Sweden&9,340,682&20&17&17&17\\
15&Austria&8,375,290&19&16&16&15\\
\seprule
16&Bulgaria&7,563,710&18&15&15&14\\
17&Denmark&5,534,738&13&12&12&12\\
18&Slovakia&5,424,925&13&12&12&12\\
19&Finland&5,351,427&13&12&12&12\\
20&Ireland&4,467,854&12&11&11&11\\
\seprule
21&Croatia&4,425,747&--&--&11&11\\
22&Lithuania&3,329,039&12&10&9&9\\
23&Latvia&2,248,374&9&8&8&8\\
24&Slovenia&2,046,976&8&8&8&8\\
25&Estonia&1,340,127&6&7&7&7\\
\seprule
26&Cyprus&803,147&6&6&6&6\\
27&Luxembourg&502,066&6&6&6&6\\
28&Malta&412,970&6&6&6&6\\
29&Iceland&317,630&--&--&--&6\\
\midrule
&{\it Total}&505,529,172&751&751 &751&754\\
\bottomrule
\end{tabular*}
\bigskip
\small\caption{The column labelled `27 States' is the Cambridge Compromise with the present 
European Union.
The next two columns include Croatia and Iceland in that order.
The divisors are 819,000 (27 States), 835,000 (28 States), 844,000 (29 States).}
\label{table2}
\end{table}

\subsection{Divisors or D'Hondt?} \label{sec:dhondt}

Democracies have extensive experience of voting systems, and a variety
of nomenclature has evolved. The following trans-Atlantic translation chart is included 
here.

\medskip
{\centering\begin{tabular*}{0.6\textwidth}{@{\extracolsep{\fill}} ccc}
\emph{rounding} & \emph{Europe} & \emph{USA}\\
\midrule %\hline
downwards & D'Hondt & Jefferson\\
standard  & Sainte-Lagu\"e & Webster\\
upwards  &      & Adams\\
\end{tabular*}

}
\medskip

The Cambridge Compromise may be reformulated as a system of any of these three types,
and we illustrate this with the case of D'Hondt's method.
Allocate to every State the minimum $m$ seats (currently $m=6$). The remaining seats are allocated
according to  D'Hondt's method subject 
to the condition that,
when any State attains a total of 96 seats, then it receives no further seats.
The ensuing allocation is identical to that of the Cambridge Compromise.

The  better to aid the reader, we give a brief explanation of the relevant  D'Hondt method in the 
presence of an integral base and maximum.
Write $B$ (\resp, $M$) for the base (\resp, maximum) allocation, and $H$ for the house-size. 
Let the population-sizes be $p_1,p_2,\dots,p_n$. 
\begin{numlist}
\item
At stage 0, allocate $B$ seats to every State.  The remaining $R=H-nB$ seats 
will be allocated sequentially as follows, until
none remain. 
\item
Suppose, at some stage, that State $i$ has been allocated $a_i$ seats in all.
Find a State $j$ such that $p_j/(a_j-B+1)$ is a maximum, and allocate the next seat to this State. 
\item 
Repeat the previous step until no seats remain, subject to the condition that any State 
achieving the maximum number $M$ of seats is removed from the process. 
\end{numlist}
It may be checked that the outcome agrees with the system
$B+$D, which was shown in Section \ref{sec:choice} to be equivalent to the
Cambridge Compromise with base $b=B-1$. Similar algorithms
are of course valid for the Sainte-Lagu\"e and Adams methods.

Ties can occur in the above algorithm, and these correspond to the non-existence of a divisor
for some house-size in the formulation of Section \ref{sec:ccbp}. There  are standard ways of breaking ties by
casting lots. However, ties are very unlikely to occur in
instances of the EU apportionment problem since populations are large and varied. 
Indeed, subject to a reasonable probabilistic
model for population-sizes, the probability of a tie may estimated rigorously. 

For further reading, see  \cite{sj11}, or perhaps \cite[p.\ 99]{baly},

\subsection{Choosing the minimum and maximum}

The better to understand the role of the minimum, CAM discussed how the minimum and base could
be reduced as further States accede to the Union.  No final recommendation was reached but two
Schemes emerged. 

In Scheme A, a cap is introduced on the proportion of seats allocated via
the minimum, and the value of the minimum is taken as large as possible subject to this cap.
For example, there are currently $27\times 6 = 162$ seats allocated thus, a proportion of
approximately 22\%. If, for example, one caps this at 25\%, the minimum remains
at 6 for a larger Union of 27--31 States, and is reduced to 5 for 32--37 States, and so on.  
The base $b$ might either be one 
fewer than the minimum (with rounding upwards), or 
might follow a rule of the type:
$b$ is the smallest fraction such that
the smallest State receives exactly the minimum number of seats (with rounding upwards, say).

In Scheme B, one determines the base as a function of the number $n$ of States, and current
practice indicates a formula of the type $b=135/n$. 
This has the advantage of decreasing steadily
as $n$ increases. However, the associated
\emph{minimum} decreases in a manner that is sensitive to the smallest
population.

Since each State receives by necessity an integral number of seats, one
effect of the allocation of seats to new States is a notable lumpiness
at the upper end of the population chart. With the minimum held constant, the
seats granted to an acceding State are taken from other States in proportion to their populations, 
and thus mostly from the larger States. 
Conversely,  any adjustment downwards in the minimum allocation
releases seats for proportional distribution between the States, of which
the largest States gain most.

CAM recommended that consideration be given to the manner 
in which the minimum allocation should vary 
in the light of changes to the European Union, and also that the functioning of
the maximum allocation be reviewed prior to future apportionments.

\subsection{Population statistics}

Census data is key to the allocation of seats in the European Parliament.
Such population data is usually collected
only once a decade. Both the year of the census and the manner of updating
can vary between countries. In addition, there can be national
variation in the definition of a resident.
CAM's final recommendation was that the European Commission be encouraged to
ensure that Eurostat
review the methods used across the Union.

\section{Summary of Recommendations}
\label{sec:recom}

\noi
\emph{Principal recommendations}
\begin{numlist}
\item Adopt the revised definition of degressive proportionality
proposed in Section \ref{sec:dp2} above.
\item For future apportionments of the European Parliament, the
method \emph{\bp} should be employed.
\item The base should be 5, and fractions should be rounded upwards.
\end{numlist}

\noi
\emph{Further recommendations}
\begin{Letlist}
\item Due consideration should be given to the manner in which the
minimum, currently 6, and base should vary in the light of future
changes in the number of Member States in the European Union.
\item The European Parliament should review the manner of
functioning of the maximum constraint on number of seats,
currently 96, prior to future apportionments.
\item The Commission should be encouraged to ensure that Eurostat
reviews the methods used by Member States in calculating their
current populations, in order to ensure accuracy and consistency.
\end{Letlist}

\section{Debate in the AFCO Committee}
\label{sec:afco}

The timetable of discussion in Brussels was as follows.
In advance of completion of the final CAM Report,
the author was invited (as Director of CAM) to deliver a preview
to the Committee on Constitutional Affairs (AFCO)  in Brussels
on 7 February 2011.  There was a Committee discussion on 15 March.  
The Rapporteur, Andrew  Duff, tabled a proposal \lq\lq for a modification of the Act concerning the election 
of the Members of the European Parliament by direct universal suffrage of 20 September 1976'',
and this was the subject of amendments by Committee members, leading in turn to
a set of so-called \lq\lq Compromise Amendments'' from the Rapporteur\footnote{Video
recordings  of the two meetings may be found at \url{http://tinyurl.com/5s63d8r}.
Versions of the proposals and amendments may be consulted at \url{http://tinyurl.com/6bzedza}.}.
A vote was taken on 19 April 2011.

The initial responses of  Committee members to the CAM recommendation
varied between curiosity verging
on support, a desire for clarification, simple misunderstanding, and downright opposition. 
Several members expressed dismay at
the \lq\lq political'' challenges of such a reorganization, and everyone was doubtless sensitive
to the needs of Member States, Political Groups, and individual Members of the European Parliament.  
Amongst the issues that stimulated some MEPs were the changes in allocations to
Member States with populations in the 7--11 million range, and the claim
by one MEP of unfair treatment of the largest Member State.

The five week intermission between the two Committee meetings permitted a period of
reflection and analysis, and contributions at the second meeting were 
generally more refined. There was some agreement 
in principle on the desirability of a formulaic approach
to apportionment, but only one speaker (apart from the Rapporteur) spoke in support of the Cambridge Compromise.
Representatives of several medium-sized countries were particularly implacable.  

Committee members tabled 138 amendments to the 
Rapporteur's Proposal for a modification of the relevant Act. The final three were
proposals to employ, \resp, the Cambridge Compromise, a parabolic method, and a power method.
These three amendments were not destined to survive the vote, presumably
as the consequences of  formulaic approaches became clearer to
some members of the Committee and of Parliament. 

Two of the Rapporteur's twelve  \lq\lq Compromise Amendments'' were 
directly relevant to the
Cambridge Compromise. Amendment B proposed a formal definition of degressive
proportionality along the lines of Section \ref{sec:dp2}, while withdrawing the proposal
to adopt a specific mathematical approach. Amendment F compressed the discussion
of a \lq\lq mathematical formula'' as follows:
\begin{quotation}
[The European Parliament] \lq\lq Proposes to enter into a dialogue with
the European Council to explore the
possibility of reaching agreement on a
\emph{durable}\footnote{Italics by the current author. Recall the \emph{three}  criteria
of Section \ref{sec:crit}; the criterion of \lq\lq impartiality" has been omitted.} 
and \emph{transparent} mathematical
formula for the apportionment of seats in
the Parliament respecting the criteria laid
down in the Treaties and the principles of
plurality between political parties and
solidarity among States.'' 
\end{quotation}
These Compromise Amendments were agreed by the Committee on 19 April 2011, and the amended
Proposal was duly carried.

It is not the purpose of this paper to speculate about the reasons for the unenthusiastic response
of the AFCO Committee to this proposal in particular, and to formulaic approaches
in general. Change can be tricky to manage and to explain to electorates, 
especially fundamental change requiring unanimity across
EU Member States and affecting the livelihoods and ambitions of
individual MEPs. The current allocations give preferential treatment to citizens of medium-sized States
at the expense of those of  larger States.   The tentacles of the Political Groups entangle the EU,
and alliances harness power and can frustrate change.    

There is also the problem of the largest State. According to the Treaty of Lisbon, 
no State shall receive more than 96 seats, whereas an uncapped allocation 
would currently give a  greater number to Germany. This feature of
Parliamentary structure is illuminated baldly by the Cambridge Compromise using
current population figures (the prominence of this cap will fade as the EU is enlarged).
  
It was argued by some MEPs that, in preferring a linear system,  
CAM had misunderstood the meaning of \lq \lq degressive proportionality".
Such critics considered that CAM should have designed
a formula to reproduce the current profile of Parliament.  Not only is this contrary
to the terms of reference received from the AFCO Committee, but also 
the author believes that mathematics is best not used as a tool to
legitimize blatantly political deals. 

The argument provides, however,  a clue
as to why formulaic approaches were 
disfavoured in the vote.
Calculations indicate that, as the number of Member States increases,  the
allocations of many formulaic systems
approach the simple linearity of the Cambridge Compromise. 
For example, with 29 States (including Croatia and Iceland) the allocations of both
the parabolic and power methods differ only very slightly from that of the Cambridge formula.
It seems that the mid-range bulge can be preserved only through \lq \lq political bartering",
and that the discussion of this paper will resurface in the future.

\section*{Acknowledgements}

The author thanks Andrew Duff for his  advice and support. He offers
his appreciation to the individual participants
at the Cambridge Apportionment Meeting, and especially to 
Friedrich Pukelsheim who shared in the organization of the event and read a draft
of this paper. 
Svante Janson kindly commented on
the representation of the Cambridge Compromise as the D'Hondt scheme of Section 
\ref{sec:dhondt}.

\bibliographystyle{amsplain}
\bibliography{ep}
\end{document}